\documentclass[a4paper,10pt]{article}

\usepackage[english]{babel} \usepackage{german, umlaute, amsmath,
  amssymb , latexsym, theorem, epic, epsfig, rotating}

\newcommand{\alphlist}{\begin{list}{(\alph{enumi})}{\usecounter{enumi}}}
\newcommand{\romanlist}{\begin{list}{(\roman{enumi})}{\usecounter{enumi}}}
\newcommand{\listend}{\end{list}}

\newcommand{\ld}{\ensuremath{,\ldots,}}
\newcommand{\ssq}{\ensuremath{\subseteq}}

\newcommand{\ra}{\ensuremath{\rightarrow}}
\newcommand{\equi}{\ensuremath{\Leftrightarrow}}
\newcommand{\follows}{\ensuremath{\Rightarrow}}

\newcommand{\T}{\ensuremath{\mathbb{T}}}
\newcommand{\N}{\ensuremath{\mathbb{N}}} 
\newcommand{\R}{\ensuremath{\mathbb{R}}}
\newcommand{\Z}{\ensuremath{\mathbb{Z}}}
\newcommand{\Q}{\ensuremath{\mathbb{Q}}}

\newcommand{\Proj}{\ensuremath{\mathbb{P}}}

\newcommand{\kreis}{\ensuremath{\mathbb{T}^{1}}}
\newcommand{\ntorus}[1][2]{\ensuremath{\mathbb{T}^{#1}}}

\newcommand{\nLim}{\ensuremath{\lim_{n\rightarrow\infty}}}

\newcommand{\qed}{{\raggedleft $\Box$ \\}}

\newcommand{\proof}{\textit{Proof:}}

\newcommand{\piin}{\ensuremath{\pi^{-1}}}

\newcommand{\incup}{\ensuremath{\bigcup_{i=1}^n}}

\newcommand{\ntel}{\ensuremath{\frac{1}{n}}}

\newcommand{\Thom}{\ensuremath{{\cal T}_{\textrm{hom}}}}

\newcommand{\Ath}{\ensuremath{A_{\theta}}}

\newcommand{\thx}{\ensuremath{(\theta,x)}}
\newcommand{\thom}{\ensuremath{\theta + \omega}}

\newcommand{\phiplus}{\ensuremath{\varphi^+}}
\newcommand{\phimin}{\ensuremath{\varphi^-}}

\newcommand{\pqphith}{\ensuremath{\displaystyle (\varphi^i_j(\theta))^{
 \scriptscriptstyle 1\leq i \leq p}_{ \scriptscriptstyle 1 \leq j \leq
q }  } }

\newcommand{\pqphiset}{\ensuremath{\{
\varphi_1^1(\theta), \ldots , \varphi^p_q(\theta) \} }}

\newcommand{\gamhat}{\ensuremath{\widehat{\gamma}}}

\newcommand{\Tth}{\ensuremath{T_{\theta}}}
\newcommand{\Tnth}{\ensuremath{T_{\theta}^{n}}}

\newcommand{\Tthx}{\ensuremath{T_{\theta}(x)}}
\newcommand{\Tnthx}{\ensuremath{T_{\theta}^{n}(x)}}

\newcommand{\That}{\ensuremath{\widehat{T}}}
\newcommand{\Thatth}{\ensuremath{\widehat{T}_\theta}}

\newcommand{\thhat}{\ensuremath{\hat{\theta}}}
\newcommand{\xhat}{\ensuremath{\widehat{x}}}
\newcommand{\yhat}{\ensuremath{\widehat{y}}}
\newcommand{\phihat}{\ensuremath{\widehat{\varphi}}}
\newcommand{\Ahat}{\ensuremath{\widehat{A}}}
\newcommand{\Uhat}{\ensuremath{\widehat{U}}}

\newcommand{\filled}[1]{\ensuremath{#1^{\mathit{fill}}}}

\theoremstyle{break}
\newtheorem{definition}{Definition}[section]
\newtheorem{thm}[definition]{Theorem}
\newtheorem{thmdef}[definition]{Theorem and Definition}
\newtheorem{lem}[definition]{Lemma}
  
\newtheorem{prop}[definition]{Proposition}
\newtheorem{questions}[definition]{Questions}

\theorembodyfont{\rmfamily}
\newtheorem{bem}[definition]{Remark}

\numberwithin{equation}{section}

\title{Towards a Classification for Quasiperiodically Forced Circle Homeomorphisms } 
\author{Tobias H. J\"ager\footnote{Friedrich-Alexander Universit\"at
    Erlangen-N\"urnberg} and Jaroslav Stark\footnote{Imperial College 
  London}}

\begin{document}

\maketitle

%
%
\begin{abstract}
  Poincar\'e's classification of the dynamics of homeomorphisms of the
  circle is one of the earliest, but still one of the most elegant,
  classification results in dynamical systems. Here we generalize this
  to quasiperiodically forced circle homeomorphisms, which have been
  the subject of considerable interest in recent years. Herman already
  showed two decades ago that a unique rotation number exists for all
  orbits in the quasiperiodically forced case. However, unlike the
  unforced case, no \textit{a priori} bounds exist for the deviations
  from the average rotation. This plays an important role in the
  attempted classification, and in fact we define a system as
  \textit{regular} if such deviations are bounded and
  \textit{irregular} otherwise. For the regular case we prove a close
  analogue of Poincar\'e's result: if the rotation number is
  rationally related to the rotation rate on the base then there
  exists an invariant strip (an appropriate analogue for fixed or
  periodic points in this context), otherwise the system is
  semi-conjugate to an irrational translation of the torus. In the
  irregular case, where neither of these two alternatives can occur,
  we show that the dynamics are always topologically transitive.
\end{abstract}


\section{Introduction}

One of the most fundamental goals of the theory of dynamical systems
is to classify systems according to their qualitative dynamical
behaviour. One of the earliest and most important of such results was
Poincar\'e's classification of homeomorphisms of the circle (e.g.
\cite{demelo/vanstrien:1993,katok/hasselblatt:1997}). Recall that for
any such map $f$ one can define the \textit{rotation number} $\rho$,
which measures the average speed of rotation of orbits around the
circle. Poincar\'e proved that this was the same for all orbits and
hence was an invariant of $f$. This leads to the following
classification:

\begin{itemize}
\item If the rotation number is rational $\rho = p/q$, then $f$ has a
  periodic orbit of period $q$.
  \item If the rotation number is irrational, then $f$ is
    \textit{semi-conjugate} to the irrational (rigid) rotation
    $R_{\rho}(\theta) = \theta+\rho$ (mod 1). Recall that this means
    that there is a continuous surjective $h$ such that $f \circ h = h
    \circ R_{\rho}$. The map $h$ is called a \textit{semi-conjugacy}
    of $f$ to $R_{\rho}$.
 \end{itemize}
 This picture in the irrational case was further completed by Denjoy
 (e.g. \cite{demelo/vanstrien:1993,katok/hasselblatt:1997}):
\begin{itemize}
\item If the rotation number of $f$ is irrational and the derivative
  of $f$ has bounded variation, then the semi-conjugacy $h$ is a
  homeomorphism. In such a case $h$ is called a {\em conjugacy} and we
  say that $f$ is \textit{conjugate} to $R_{\rho}$.
\end{itemize}
Subsequent work, which formed a key part of the so called KAM theory,
then characterized the conditions under which the conjugacy could be
guaranteed to be smooth or even analytic (e.g.
\cite{demelo/vanstrien:1993,katok/hasselblatt:1997}).  It is natural
to attempt to generalize these results to higher dimensions, that is
to homeomorphisms of the torus \ntorus[n]. However, even in the case
$n=2$ this turns out to be difficult.  Firstly one finds that, in
general, the rotation vector depends on the orbit. Hence, instead of a
unique, well defined rotation number one obtains a rotation set (see
\cite{misiurewicz/ziemian:1989}). Furthermore, even if this is reduced
to a single rotation vector, examples of Furstenberg (see
\cite{mane:1987}) and Herman (\cite{herman:1986,herman:1983}) show
that the dynamics may still not be compatible with either of the two
cases in Poincar\'e's classification.

Here we study a class of systems which are somewhat intermediate
between dimension one and two, namely quasiperiodically forced
systems, which have recently been the subject of considerable
interest. Specifically, we shall consider \textit{quasiperiodically
  forced circle maps}, which are continuous maps of the form
\begin{equation}                               
  \label{generalsystem}
  T : \ntorus \ra \ntorus \ , \ \thx \mapsto (\thom,\Tthx) \ ,
\end{equation}
where the \textit{fibre maps} \Tth\ (defined by $\Tthx = \pi_2 \circ
T\thx$) are all orientation-preserving circle homeomorphisms and
$\omega$ is irrational. To guarantee any required lifting properties
we additionally assume that $T$ is homotopic to the identity on
\ntorus[n] and denote the class of such systems by \Thom.

On the one hand, the skew product structure implies that such maps
have similar order-preserving properties to homeomorphisms of the
circle. This is sufficient to ensure that the rotation vector is
uniquely defined (see Thm.\ \ref{def:rotnum} below). On the other
hand, the examples in \cite{furstenberg:1961} and in
\cite{herman:1983} show that the dynamical phenomena which can be
found in \Thom \ are already much richer than those in the class of
circle homeomorphisms.

In addition, quasiperiodically forced systems are also interesting in
their own right, as they occur in various situations in physics.  For example,
the quasiperiodically forced Arnold circle map serves as a simplified model
for an oscillator forced at two incommensurate frequencies (e.g.\
\cite{ding/grebogi/ott:1989}), and the spectral theory of certain discrete
Schr\"odinger operators is intimately related to the dynamical properties of
quasiperiodically forced Moebius transformations (e.g.\ \cite{avila/krikorian}).

Our aim here is thus to derive a classification for maps of the form
(\ref{generalsystem}) which is analogous to that for
homeomorhisms of the circle. It turns out, however, that a key boundedness
property in the one dimensional case no longer holds for
quasiperiodically forced maps. In particular if $\hat{f}$ is the lift  to \R
of a homeomorhism of the circle, then the monotonicity of
$\hat{f}$ and periodicity of $\hat{f} - Id$ implies that

\begin{equation}                               
  \label{pairbound}
  |\hat{f}^n(x) - \hat{f}^n(y)| \leq |x - y| + 1 \ \forall n \in \N,\ x,y \in \R.
  \end{equation}
  This gives a simple proof of the existence of the rotation number
  $\rho$ and the following uniform bound for \textit{deviations from
    uniform rotation}
\begin{equation}                               
  \label{uniformbound}
  | \hat{f}^n(x)-x-n\rho | \leq 1 \ \forall n \in \N,\ x \in \R.
  \end{equation}
  Interestingly this holds even if $\hat{f}$ is merely monotone
  increasing, but not necessarily continuous \cite{rhodes86,rhodes91}.
  In the case of quasiperiodically forced maps, any given fibre map
  \Tth \ is a homeomorhism of the circle, and hence the analogue of
  (\ref{pairbound}) holds within each fibre. This is sufficient to
  prove the existence of a rotation number and fairly elementary
  arguments from ergodic theory then show that this rotation number is
  the same for all fibres
  \cite{herman:1983,stark/feudel/glendinning/pikovsky:2002}.
  Unfortunately, although all fibres therefore rotate at the same
  average rate, the deviations between fibres need not be bounded and
  hence the analogue of (\ref{uniformbound}) no longer need hold, even
  within a single fibre. Surprisingly, if the analogue of
  (\ref{uniformbound}) does hold for even a single orbit, then
  \cite{stark/feudel/glendinning/pikovsky:2002} show that it holds for
  all orbits. Here we call such a system \textit{regular}. The main aim of this paper is to show that for regular systems we have a classification similiar to that for circle maps. The details of this are given in Section \ref{Regularcase}
  below. We then go on to discuss the \textit{irregular} case in
  Section \ref{Irregularcase}. Whilst we are unable to give an
  elegant classification according to the rotation number, as in the
  regular case, we are at least able to show that all systems in this
  class are topologically transitive.
  
  We begin by collecting a variety of preliminary facts and
  constructions, which are required in the succeeding sections. We also
  discuss a result of Furstenberg which can be seen as an
  analogue to Poincar\'e's classification, but in a measure-theoretic
  rather than a topological sense.


\section{Preliminaries}
\label{Prerequisites}

\subsection{Notation}

In the following $m$ will denote Lebesgue measure on \kreis, $\lambda$
Lebesgue measure on \ntorus\ and $\pi_{i} : \ntorus \ra \kreis \ 
(i=1,2)$ the projection to the respective coordinate. We shall denote
the $q-$fold cover $\R / q\Z$ of the circle \kreis\ by $\kreis_q$.  If
there is no possible ambiguity, all projections from either \R\ or
$\kreis_q$ to $\kreis$ or from $\R^2$, $\kreis \times \R$, $\T^1_q
\times \T^1_k$ etc.\ to \ntorus\ will be denoted by $\pi$. Lifts of
$T$ to the covering spaces $\kreis \times \R$ or $\R^2$ will be
denoted by $\That$ and variables in the respective covers will be 
written $\thhat,\xhat$ and so on. It will also be useful to define
$\ntorus[q]_* := \{ (\theta_1 \ld \theta_q) \mid \theta_j \in \kreis,$ 
$\theta_i \neq \theta_j \ \forall i\neq j \in \{1 \ld q \} \}$.

For any set $A \subseteq \ntorus, \ \kreis \times \R$ or $\R^2$ let
$A_\theta := A \cap \pi^{-1}\{\theta\}$ be the restriction of $A$ to
the fibre above $\theta$. Thus, either $A_\theta = \{ x \in \kreis
\mid \thx \in A \}$ or $A_\theta = \{ x \in \R \mid \thx \in A \}$.
For any pair of functions $\varphi, \psi : \kreis \ra \R$ we let
$[\varphi,\psi]$ be the region $\{ (\theta,\widehat{x}) \in \kreis
\times \R \mid \varphi(\theta) \leq \xhat \leq \psi(\theta) \}$
between the graphs of $\varphi$ and $\psi$, and similarly for
functions $\kreis \ra \kreis$ or $\R \ra \R$. The notation $\varphi <
\psi$ will mean $\varphi(\theta) < \psi(\theta) \ \forall \theta \in
\kreis$ and similarly for $\varphi \le \psi, \varphi = \psi, \varphi >
\psi$ etc.

Finally, when considering fibre maps of iterates of $T$ or their
inverses we use the convention $\Tnth := (T^n)_{\theta}\ \forall
n\in\Z$, so that $\Tnthx = \pi_2 \circ T^n\thx$.


\subsection{Minimal Sets for Monotone Quasiperiodically Forced Maps}
\label{Minimalsets}


We will mainly work with the lift $\That$ of a map $T \in \Thom$
rather than with the original map $T$, so we first collect some basic
statements about such lifts. The main properties we need are the skew
product structure and the monotonicity of the fibre maps. The fact
that the lift (or more precisely $\That - (0,Id)$) is also periodic
(in $x$) will only be used in subsequent sections. Therefore,
throughout this section we assume that $\ \That : \kreis \times \R \ra
\kreis \times \R$ has skew product structure as in
(\ref{generalsystem}) and that all fibre maps $\That_\theta$ are
(non-strictly) monotonically increasing, that is $\That_\theta(x) \le
\That_\theta(y)$ for all $x \le y$ .

Obviously there cannot be any fixed of periodic points for $\That$,
and due to the minimality of the forcing irrational rotation any
compact invariant set must project down to the whole circle. Therefore
it is natural to replace fixed points by invariant graphs:

\begin{definition}[Invariant Graphs]
  \label{def:invgraphs}
  Let $\That$ be as above and suppose a function $\varphi : \kreis \ra
  \R$ satisfies
\begin{equation}
  \label{eq:invgraphs}
   \That_\theta(\varphi(\theta)) \ = \ \varphi(\thom) \ .
\end{equation}
Then $\varphi$ is called a $\That$-invariant graph.
\end{definition}
Note that in the unforced case where $\That_\theta$ is
independent of $\theta$, a fixed point of $\That_\theta$ becomes a
horizontal invariant graph for \That.  As long as no ambiguities can
arise, the point set $\Phi := \{ (\theta,\varphi(\theta)) \mid \theta
\in \kreis \}$ will also be called an invariant graph. There is a
natural relation between compact invariant sets and invariant graphs
(see \cite{stark:2003}):
\begin{lem}
  \label{lem:boundinggraphs}
  Let $A \subset \kreis \times \R$ be a compact $\That$-invariant set.
  Then
\[
    \phiplus_A(\theta) := \sup \Ath \ \ \textrm{ and
    } \ \ \phimin_A(\theta) :=  \inf \Ath
\]
both define invariant graphs. Furthermore $\phiplus_A$ is upper
semi-continuous and $\phimin_A$ is lower semi-continuous.
\end{lem}

\ \\ If we apply this procedure to an invariant graph $\varphi$ with
compact closure $\overline{\Phi}$, then we simplify notation by
writing $\phiplus := \phiplus_{\overline{\Phi}}$ and $\phimin :=
\phiplus_{\overline{\Phi}}$. Further, we write $\varphi^{+-}$ instead
of $(\varphi^+)^-$, etc. Particularly interesting is the case where
$A$ is a minimal set in the sense of topological dynamics (so that
every orbit is dense in $A$, e.g. \cite{katok/hasselblatt:1997}); see
\cite{stark:2003} again:

\begin{lem}
  \label{lem:minimalsets}
Let $A$ be a minimal set for $\That$. Then $\varphi_A^{+-} = \phimin_A, \ 
\varphi_A^{-+} = \phiplus_A$ and $\overline{\Phi^+_A} = \overline{\Phi^-_A}
= A$. Furthermore, the set $\{ \theta \in \kreis \mid \phiplus_A(\theta) =
\phimin_A(\theta) \}$ is residual.
\end{lem}
This motivates the following definition:
\begin{definition}[Pinched Sets, Reflexive Graphs]
  \alphlist
  \item A compact set $A \subseteq \kreis \times \R$ with $\pi_1(A) = \kreis$
   is called pinched if $\{ \theta \in \kreis \mid \phiplus_A(\theta) =
   \phimin_A(\theta) \}$ is residual.
  \item An upper semi-continuous graph $\varphi : \kreis \ra \R$ is called
   reflexive if $\varphi^{-+} = \varphi$. Similarly, a lower semi-continuous
   graph $\varphi$ is called reflexive if $\varphi^{+-} = \varphi$.
 \listend
\end{definition}
We can thus restate Lemma \ref{lem:minimalsets} in the following way: a
compact $\That$-invariant set $A$ is minimal if and only if it is the closure
of a reflexive graph, and the bounding graphs of minimal sets are always
reflexive. Furthermore, every minimal set is pinched. 

There is a simple way of producing reflexive graphs from
semi-continuous ones:
\begin{lem}\label{lem:reflexivegraphs}
  Suppose $\varphi : \kreis \ra \R$ is upper semi-continuous. Then $\varphi^-$
  is reflexive. Similarly, if $\varphi$ is lower semi-continuous then
  $\varphi^+$ is reflexive. 
\end{lem}
We prove the case where $\varphi$ is upper semi-continuous, the
lower-semicontinuous case is analogous. We have to show that
$\varphi^- = \varphi^{-+-}$. First of all, as $\varphi^- \leq \varphi$
and due to the semi-continuity of the two graphs, the set
$[\varphi^-,\varphi]$ is compact and contains $\Phi^-$. Therefore it
also contains $\Phi^{-+}$, so that $\varphi^- \leq \varphi^{-+} \leq
\varphi$ and similarly $\varphi^- \leq \varphi^{-+-} \leq \varphi$.
But, as $\varphi^{-+-}$ is lower semi-continuous, the set
$[\varphi^{-+-},\varphi]$ is compact as well, consequently it contains
$\overline{\Phi}$ and thus also $\Phi^-$. This proves the reverse
inequality $\varphi^{-+-} \leq \varphi^-$.

\qed

\ \\ Another nice property of minimal sets is that they are strictly
ordered. We introduce the following notation to describe the order of
sets: If $A,B$ are bounded subsets of $\kreis \times \R$ with
$\pi_1(A) = \pi_1(B) = \kreis$, then
\begin{eqnarray*}
  A \ \preccurlyeq \ B & :\equi & \phimin_A \ \leq \ \phimin_B \ \
  \textrm{and} \ \ \phiplus_A \ \leq \ \phiplus_B \\
  A \ \prec \ B & :\equi & \phiplus_A < \phimin_B \ .
\end{eqnarray*}
We then have (see \cite{stark:2003} again):
\begin{lem}
  \label{lem:ordering}
If $A,B \subseteq$ are minimal w.r.t.\ $\That$, then either $A\prec B$, $A=B$
or $A \succ B$.
\end{lem}
Using Lemma \ref{lem:reflexivegraphs}, it is now easy to see that for any compact
invariant set $A$ the `highest minimal set' it contains is given by
$\overline{\Phi_A^{+-}}$.

Note that if $A$ is a compact invariant set then so is
$[\phimin_A,\phiplus_A]$ but the two sets are not neccessarily equal.
As it will sometimes be convenient to work with the second,
\textit{filled-in} set rather than with the first, we make the
following definition:
\begin{definition}[Filled-In Sets, Minimal Strips]
  If $A \subseteq \kreis \times \R$ is a compact set with $\pi_1(A) = \kreis$,
  we denote the filled-in set by 
 \[
     A^{\mathit{fill}} := [\phimin_A,\phiplus_A] \ .
 \]
 $A$ is called a strip if it consists of one interval on every fibre
 (i.e.\ $A = \filled{A}$). If the two bounding graphs are reflexive,
 $\filled{A}$ is called a minimal strip.
\end{definition}

\begin{bem}
  \label{bem:rtwosets}
\alphlist
\item As there can be no other semi-continuous graph between a pair of
  reflexive graphs, a minimal strip cannot contain any smaller strip.
\item Note that \textit{a priori} minimality of strips is a purely
  topological property and not neccessarily related to any dynamics.
\item Similarly to minimal sets, the intersection of two invariant
  minimal strips is either empty or the two strips are equal. However,
  in the case of strips we do not even have to invoke the dynamics:
  two minimal strips are already equal if their intersection projects
  down to the whole circle. This follows from the fact that the
  intersection is a strip and a minimal strip cannot strictly contain
  any smaller strip.
\item As minimal sets, two disjoint strips are strictly ordered. This
  can be seen as follows: If $A,B$ are both strips, then due to the
  semi-continuity of the bounding graphs both $\{\theta \in \kreis
  \mid \phiplus_A < \phimin:B \}$ and $\{ \theta \in \kreis \mid
  \phiplus_B < \phimin_A \}$ are open. As both sets together cover the
  whole circle, one of them must be empty.
\item Below we will also work with lifts to $\R^2$ instead of $\kreis
  \times \R$. We call a set $A \subseteq \R^2$ a strip if it consists
  of one interval on every fibre, projects down to all of $\R$ and
  $\piin_1(K) \cap A$ is compact for every compact set $K \subseteq
  \R$. The above definitions of bounding graphs, reflexivity,
  filled-in sets, minimal strips and the order relations between sets
  can all be applied to strips in $\R^2$ in the same way.  \listend
\end{bem}


\subsection{Invariant Strips and Regular Invariant Graphs}
\label{Invariantstrips}

Invariant graphs are the natural analogue of fixed points for
quasiperiodically forced monotone maps. When we come to study
quasiperiodically forced circle maps we also need to generalize the
concept of periodic orbits to that of periodic graphs. Since such
graphs may wrap around the torus more than once in the
$\theta$-direction, the following definition is a little bit more
complicated than Def.\ \ref{def:invgraphs}.

\begin{definition}[$p,q$-Invariant Graphs]
  \label{def:invgraph}
  Let $T \in \Thom$, $p,q \in \N$. A $p,q$-invariant graph is a
  measurable function $\varphi : \kreis \ra \ntorus[pq]_*\ , \ \theta
  \mapsto \pqphith$ with
  \begin{equation}
    \label{eq:invgraph}
    \Tth(\pqphiset ) \ = 
    \ \{ \varphi_1^1(\thom) , \ldots , \varphi^p_q(\thom) \} 
  \end{equation}
  for $m$-a.e. $\theta \in \kreis$, which satisfies the following:
  \romanlist
  \item $\varphi$ cannot be decomposed (in a measurable way) into
  disjoint subgraphs $\varphi^1 \ld \varphi^m \ (m\in\N)$ which
  also satisfy (\ref{eq:invgraph}).
  \item $\varphi$ can be decomposed into $p$ $p$-periodic $q$-valued
    subgraphs $\varphi^1,\ldots,\varphi^p$.
  \item The subgraphs $\varphi^1 \ld \varphi^p$ cannot further be
    decomposed into invariant or periodic subgraphs.
  \listend
  If $p=q=1$, $\varphi$ is called a simple invariant graph.  The point
  set $\Phi := \{ (\theta,\varphi^i_j(\theta)) \mid \theta \in \kreis,\
  1\leq i \leq p,\ 1 \leq j \leq q \}$ will also be called an
  invariant graph, but labeled with the corresponding capital letter.
  An invariant graph is called continuous if it is continuous as a
  function $\kreis \ra \ntorus[pq]_*$.
\end{definition}
For a further discussion and some elementary examples see
\cite{jaeger/keller:2003}. The above definition only assumes
measurability, and thus an invariant graphs does not neccessarily have
any topological structure. This leads to the difficulty that it is
not possible to define a lifts of a  $p,q$-invariant graph in any obvious way.
Despite this, the above concept leads to a measure-theoretic
classification which was already given by H.\ Furstenberg in 1961 
\cite{furstenberg:1961}, albeit using very different terminology:
\begin{thm}[Furstenberg, \cite{furstenberg:1961} Thm.\ 4.1]\label{thm:furstenberg}
  Let $T \in \Thom$. Then either there exists a $p,q$-invariant graph
  and every invariant ergodic measure is supported on such a graph, or
  $T$ is uniquely ergodic and isomorphic to a uniquely ergodic skew
  translation of the torus, i.e.\ a map of the form $A: \thx \mapsto
  (\thom,x+a(\theta))$.
\end{thm}
Recall that two dynamical systems $(X,T)$ and $(Y,S)$ with invariant
probability measures $\mu$ and $\nu$ are called isomorphic, if there
exist sets $A\subseteq X$ and $B \subseteq Y$ with $\mu(A)=\nu(B)=1$
and a bijective measurable map $h:A \ra B$ which maps $\mu$ to $\nu$
and satisfies $h \circ T = S \circ h$. Although not explicitly
mentioned in \cite{furstenberg:1961}, Furstenberg's construction gives
a lot of additional information about the isomorphism $h$ in the
situation of Thm.\ \ref{thm:furstenberg}. This will be discussed in
some detail at the end of Section \ref{Irregularcase}.

Note that this result still holds if the irrational rotation on the
base is replaced by any uniquely ergodic dynamical system. The
weakness of this classification is that it does not for example
distinguish between a continuous invariant curve and an invariant
graph which is dense in the whole torus. These are clearly different
situations from a topological point of view.  On the other hand, it
would be far too restrictive to consider only continuous invariant
graphs since numerical experiments suggest that it is possible for
instance to have regular systems with a rational rotation number, but
no continuous invariant graph \cite{feudel/kurths/pikovsky:1995}. In
such cases the invariant graph appears to have been replaced by a more
complex set, which often has bifurcated from a continuous invariant
graph. We have already seen a good candidate for such a set in the
previous section in the form of an invariant strip, with its related
semi-continuous invariant bounding graphs. This is exacly the concept
we want to use now in order to complement Furstenberg's result with a
topological classification. Note that there is no straightforward way
to define semi-continuity for a function $\varphi : \kreis \ra
\kreis$. In order to obtain some kind of semi-continuity property for
the bounding graphs of a strip in \ntorus\ we therefore include the
existence of an appropriate `reference curve' in the definition of an
invariant strip. Such a curve also ensures the lifting properties to
the cover that we shall require. An appropriate definition is given by
\begin{definition}[$q$-Curves]
  \label{def:qcurves}
Let $\gamhat : \R \ra \R$ be a continuous function which satisfies
\begin{equation}
  \label{eq:qcurve}
     \gamhat(\thhat+q) \ = \ \gamhat(\theta) + k \ \ \textrm{ and } \ \
     \gamhat(\thhat + l) - \gamhat(\thhat) \notin \Z \ \ \forall \thhat \in
     \R, \ l = 1 \ld q - 1    
\end{equation}
for some $q \in \N$ and $k \in \Z$. A $q$-curve is then the projection 
\[
    \gamma : \kreis \ra \ntorus[q]_*, \ \theta \mapsto (\gamma_1(\theta) \ld \gamma_q(\theta))
\]
with $\gamma_i(\theta) := \pi(\gamhat(\thhat + i - 1))$. The point set
$\Gamma := \{ (\pi(\theta),\gamma(\theta)) \mid \theta \in [0,q), \ 0 \leq i < q \} 
\subseteq \ntorus$
will equally be refered to as a $q$-curve if no ambiguity results. If $q=1$, $\gamma$ will be called a simple curve.
\end{definition}
Such a $q$-curve $\gamma$ has the property that it is continuous as a function
from \kreis\ to $\T^q_*$ and consists of one connected component only (as a
point set).  Conversely, it follows from the usual arguments for the existence
of lifts that any graph $\gamma : \kreis \ra \ntorus[q]_*$ with these
properties is a $q$-curve, as it can be lifted to a function $\gamhat : \R \ra
\R$ which satisfies (\ref{eq:qcurve}) and projects down to $\gamma$.

\begin{definition}[$p,q$-Invariant Open Tubes, $p,q$-Invariant Strips]
  \label{def:tubes}
  Let $T \in \Thom$.
\alphlist
\item An open $T$-invariant set $U$ is called a $1,q$-invariant open
  tube if there exists a $q$-curve $\gamma$ such that for all $\theta
  \in \kreis$, $U_\theta$ consists of $q$ disjoint open and nonempty
  intervals $U_{\theta,j}$ such that $\gamma_j(\theta) \in
  U_{\theta,j} \ \forall j = 1 \ld q$.
\item When $U^1$ is a $1,q$-invariant open tube for $T^p$ containing
  the $q$-curve $\Gamma^1$ and the sets $U^i := T^{i-1}(U^1) \ (i=2\ld
  p)$ are pairwise disjoint, then $U:= \bigcup_{i=1}^{p} U^i$ is
  called a $p,q$-invariant open tube, each of its connected components
  $U^i$ containing the $q$-curve $\Gamma^i := T^{i-1}(\Gamma^1)$.
\item A compact invariant set $A$ is called a $p,q$-invariant strip if
  its complement is a $p,q$-invariant open tube. If $p=q=1$ then $A$ is called a 
simple invariant strip.
\item $\varphi : \kreis \ra \ntorus[pq]_*$ is called a regular
  $p,q$-invariant graph if it is the bounding graph of a
  $p,q$-invariant open tube, i.e.\ 
\[
\varphi^i_j(\theta) = \gamma^i_j(\theta) + \inf \{x \in [0,1) \mid
\gamma^i_j(\theta) + x \in U_\theta^c \}
\]
or
\[
\varphi^i_j(\theta) = \gamma^i_j(\theta) - \inf \{ x \in [0,1) \mid
\gamma^i_j(\theta) - x \in U_\theta^c \} \ ,
\]
where $U = \bigcup_{i=1}^{p} U^i$ and the $\gamma^i$ are $q$-curves
contained in $U^i$. If $p=q=1$ then $\varphi$ is called a simple regular
invariant graph.
\item A $p,q$-invariant graph is called irregular, if it is not regular.
\listend
\end{definition}
Of course, a continuous invariant graph is a special case of a regular
invariant graph and also of an invariant strip.

As mentioned already, using the $q$-curves contained in an open tube we can
lift these objects to the covering spaces $\kreis \times \R$ or
$\R^2$: First of all, if $\gamma$ is a $q$-curve then there are $q$
different lifts $\gamhat^1 < \ldots < \gamhat^q$ such that
$\gamhat^i(0) \in [0,1)$. From these we can produce all possible lifts
by the addition of an integer, i.e.\ $\gamhat^{i+nq} := \gamhat^i +
n$. If $\gamma$ belongs to the $1,q$-invariant tube $U$ we can lift
$U$ `around' the respective lifts of $\gamma$ and obtain, up to
addition of integers, $q$ different lifts $\widehat{U}^1 \ld
\widehat{U}^q$ of $U$. The important thing now is to see that these
also have a certain invariance property: they might not be invariant
under a lift of $T$ itself, but at least it is possible to choose a
lift of $T^q$ which leaves them invariant. Of course, the same holds
for $1,q$-invariant strips.

In order to see this we choose a lift $\That$ of $T$ such that
$\That(0,\gamhat_1(0)) \in \widehat{U}^i$ for some $i \in \{ 1 \ld q \}$ and claim
that $\That(\hat{\Gamma}^1) \subseteq \widehat{U}^i$, which then immediately implies
$\That(\widehat{U}^1) = \widehat{U}^i$. Otherwise let $\thhat^* := \inf\{ \thhat \geq
0 \mid \That({\thhat},\gamhat^1(\thhat)) \notin \widehat{U}^i \}$. Then
$\That(\thhat^*,\gamhat^1(\thhat^*)$ is contained in some other lift
$\widehat{U}^*$ of $U$. But as $\widehat{U}^*$ is open this means that the set contains
$\That(\thhat,\gamhat^1(\thhat))$ for all $\thhat$ from a whole neighbourhood
of $\thhat^*$, contradicting the definition of $\thhat^*$.
Now the monotonicity of the fibre maps implies $\That(\widehat{U}^n) =
\widehat{U}^{n+i-1} \ \forall n \in \N$ and consequently
$\That^q(\widehat{U}^1) = \widehat{U}^1 + (0,i)$. But this means we can choose
another lift of $T^q$ such that the $\widehat{U}^n$ are left invariant.
 
Now we can prove two elementary statements about invariant strips,
which will simplify dealing with these objects later on.

\begin{lem} \label{lem:simplegraph}
Let $T\in \Thom$. Then every compact invariant set $A$ which consists of exactly
one non-trivial (i.e.\ $A_\theta \neq \emptyset$ or \kreis) interval on every fibre is
a simple invariant strip.
\end{lem}
\proof \\
Obviously $U := A^c$ is an open set which also consists of exactly one non-trivial
interval on every fibre. It remains to show that there exists a simple curve in
$U$.

To that end, note that as $U$ is open and $U_\theta \neq \emptyset$ for any
$\theta \in \kreis$ we can choose an open box $W_\theta \subseteq U$ with
$\theta \in I_\theta := \pi_1(W_\theta)$. By compactness, finitely many
intervals $I_{\theta_1} \ld I_{\theta_n}$ cover the whole circle. Let $I_j'
\subseteq I_{\theta_j}$ be disjoint open intervals with $\kreis \subseteq
\overline{\incup I_j'}$. Over any of these intervals there is a constant line
segment contained in $U$, and as $A$ only consists of one inverval on every
fibre we can always join two adjacent segments by a vertical line in one of
the two directions on the circle. That way we get a closed line in $U$
which wraps once around the torus in the $\theta$-direction. Due to the
vertical parts this line is not yet the graph of a continuous function, but as $U$ is
open we can slightly tilt them without leaving $U$  to obtain the simple curve
we require.

\qed

\begin{lem} \label{lem:liftedgraph}
$T \in \Thom$ has a $p,q$-invariant strip if and only if there exist numbers
$n,l,k \in \N$ such that a lift of $T^n$ to $\kreis_l \times \kreis_k$ has a
simple invariant strip.
\end{lem}
\proof \\ {``$\follows$'':} As $T$ has a regular $p,q$-invariant
graph there exists a $1,q$-invariant open tube $U$ for $T^p$. As argued above,
a lift $\widehat{U}$ of $U$ to $\R^2$ is invariant under a suitable lift
of $T^{pq}$. If $\widehat{\Gamma} \subseteq \widehat{U}^1$ is a lift of the
corresponding $q$-curve and $\gamhat(\thhat + q) = \gamhat(\thhat) + k$, then
$\widehat{U}_{\thhat + q} = \widehat{U}_{\thhat} + k$. Consequently the
projection of $U$ onto $\kreis_q \times \kreis_k$ is a simple invariant open
tube with respect to a suitable lift of $T^{pq}$.

\ \\ {``$\Leftarrow$'':} Let $F : \kreis_l \times \kreis_k$ be the
lift of $T^n$ and $\tilde{A}$ the corresponding simple invariant strip.  W.l.o.g.\ we
can assume that $\tilde{A}$ does not contain any other invariant strip. (Otherwise we
lift $\tilde{A}$ to $\R^2$, take the highest minimal strip it contains as described in Section
\ref{Minimalsets} and project it down again. This new strip will have the
required property then.) As $\tilde{A} + (0,1)$ is a minimal invariant strip as well,
it is disjoint from $\tilde{A}$ (unless $k = 1$). Hence the projection of $\tilde{A}$ to
$\kreis_l \times \kreis$ consists of one non-trivial interval on every fibre,
and using Lemma \ref{lem:simplegraph} we see that it is a simple invariant strip for
the corresponding lift of $T^n$. Therefore we can assume $k = 1$. Likewise we
can assume that $l$ is minimal, i.e.\ that $\tilde{A} + (i,0)$ is disjoint
from $\tilde{A}$ for all $i= 1 \ld l-1$, otherwise the two strips are equal
and we can project down to $\kreis_i \times \kreis$. 

We now need to show that the projection $A :=
\pi(\tilde{A})$ of $\tilde{A}$ to $\ntorus$ is a $1,l$-invariant strip
for $T^n$. The problem is that we cannot simply project the $q$-curve
in $\kreis_l \times \kreis$ down, as it might then intersect $A$.
Therefore we have to construct a new $q$-curve in the complement of
$A$. In order to do so we choose a lift $\widehat{A}$ of $\tilde{A}$
in $\R^2$. As $\Uhat := (\phiplus_{\Ahat},\phimin_{\Ahat}+1)$ is a
fundamental domain of $\tilde{A}^c$, there exists a unique lift of
every set $\tilde{A} + (i,0)$ in $\Uhat$. All these lifts are minimal
strips, and hence they are strictly ordered (see Remark 
\ref{bem:rtwosets}(b)). Furthermore, their union equals $\piin(A) \cap
\widehat{U}$. Obviously there exists a continuous curve $\gamhat$
between $\Ahat$ and the first of these strips above $\Ahat$, and we can
choose \gamhat\ such that it is the lift of a $q$-curve in the
complement of $A$.

\qed


\subsection{Fibrewise Rotation Numbers and $\rho$-Bounded Orbits} \label{Rotationnumbers}
As already indicated in the introduction, the monotonicity of each
fibre map combined with the unique ergodicity of the forcing rotation
ensure the existence of a fibrewise rotation number:

\begin{thmdef}[Herman, \cite{herman:1983}]
  \label{def:rotnum}
  Let $T \in \Thom$  and $\That : \kreis \times  \R \ra  \kreis \times
  \R$ be a lift of $T$. Then the limit 
  \begin{equation}
    \label{eq:rotnum}
        \rho_{\That} := \nLim \ntel(\Thatth^n(\xhat) - \xhat) \ 
  \end{equation}
  exists and is independent of $\theta$ and \xhat, the convergence in
 (\ref{eq:rotnum}) is uniform on $\kreis \times \R$ and in addition
  \begin{equation}
  \label{eq:rotnumII}
  \rho_{\That} = \nLim \ntel \int_{\kreis} \widehat{T}_\theta^n(0) \
  d\theta  
  \end{equation}
Further more, $\rho_T := \rho_{\That} \ \bmod 1$ is independent of the choice
 of the the lift \That.  It is called the fibrewise rotation number of $T$.
\end{thmdef}
However, there is a crucial difference between the quasiperiodically
forced case and an unforced homeomorphism of the circle. As
highlighted in the introduction, in the latter case there is a bound
(\ref{uniformbound}) on the possible deviations of any orbit from the
average rotation. For a general quasiperiodically forced map such a
bound need not exist, motivating the following definition:
\begin{definition}[$\rho$-Bounded Orbits]  \label{def:boundedorbits}
\alphlist
\item If $\That$ is a lift of $T \in \Thom$ and $\rho \in \R$, we say that the
      orbit of $(\theta,\xhat)$ is $\rho$-bounded if there exists a constant
      $C > 0$ such that $|\That_\theta^n(\xhat) - \xhat - n\rho| \leq C \
      \forall n \in \N$.
\item Similarly, we say that the orbit of $(\theta,\xhat)$ is $\rho$-bounded
      above (below) if there exists $C > 0$ such that $\That_\theta^n(\xhat) -
      \xhat - n\rho \leq C \ \ (\geq -C) \ \forall n \in \N$.
\item If the orbits of $\That$ are $\rho_{\That}$-bounded, then the same is
      true for any other lift of $T$. Thus we say $T$ has $\rho_T$-bounded
      orbits if the lifts have bounded orbits in the above sense.
\listend
\end{definition}
More informally we will also speak of the \textit{deviations from the constant
 rotation} when refering to the quantities $|\That_\theta^n(\xhat) - \xhat -
 n\rho_{\That}|$. The following result about boundedness of orbits is taken from
 \cite{stark/feudel/glendinning/pikovsky:2002}:
\begin{thm}
  \label{thm:boundedness}
Let $T \in \Thom$. If there exists one $\rho_T$-bounded orbit, then all orbits
are $\rho_T$-bounded and the constant $C$ in (\ref{def:boundedorbits}) can be
chosen uniformly for all $\thx \in \ntorus$. If there is no $\rho_T$-bounded
orbit, there exists at least one orbit which is $\rho_T$-bounded above and one
which is $\rho_T$-bounded below, but for a residual set of $\theta$'s all
orbits on the $\theta$-fibre $\{\theta\} \times \kreis$ are both
$\rho_T$-unbounded above and below.
\end{thm}
Note that the uniformity of $C$ is neither explicitly stated nor proved in  \cite{stark/feudel/glendinning/pikovsky:2002}. However, Lemma 7 of this reference does show that there is a uniform constant for all $\thx \in U \times \kreis$ where $U$ is a residual set, and extending this to the whole of \ntorus\ is elementary.
As we shall see in the next section, the analogue of Poincar\'e's classification for circle maps only applies to systems where all orbits are $\rho$-bounded. It is thus reasonable to define

\begin{definition}[Regular Systems]
A map $T \in \Thom$ is called regular if all orbits are
$\rho_T$-bounded. Otherwise it is called irregular.
\end{definition}
Finally, we relate the existence of a $p,q$-invariant strip to the properties of the rotation number.

\begin{definition}[Rational Dependence]
  \label{def:rationaldep}
Two numbers $\omega,\rho \in \R$ (or \kreis) are said to be rationally dependent 
if there exists 
$(l,k,q) \in \Z^3 \setminus \{(0,0,0)\}$ such that $l+ k\omega + q\rho
= 0$. Otherwise they are called rationally independent.  
\end{definition}

\begin{definition}[Irrational Torus Translation]
For $\omega,\rho \in \kreis$ let $R_{\omega,\rho} : \thx \mapsto (\thom,
x+\rho)$. We call $R_{\omega,\rho}$ an irrational torus translation if
$\omega$ and $\rho$ are not rationally dependent. 
\end{definition}
Note that as we assumed our rotation number $\omega$ to be irrational,
rational dependence of $\omega$ and $\rho_T$ is equivalent to the
existence of $l,k,q \in \Z, \ q \neq 0$ such that $\rho_T =
\frac{k}{q}\omega + \frac{l}{q}$. The following result, which relates
the concepts of invariant strips and fibrewise rotation numbers, can
be found in \cite{jaeger/keller:2003} (Lemma 3.9).

\begin{prop}
  \label{lem:tubes}
Let $T \in \Thom$ and suppose there exists a $p,q$-invariant strip. Then
$\omega$ and $\rho_T$ are rationally dependent. Furthermore, the orbits of $T$
are $\rho_T$-bounded.
\end{prop}


\section{The Regular Case}
\label{Regularcase}

For regular systems the
following statement is in perfect analogy with Poincar\'e's classical result:

\begin{thm}\label{thm:poincare}
Suppose $T \in \Thom$ is regular. Then one of the following holds:
\alphlist
\item $\rho_T$ and $\omega$ are rationally independent and there exists a
  regular invariant graph for $T$.
\item $\rho_T$ and $\omega$ are rationally dependent and $T$ is
 semi-conjugate to the irrational torus-translation
 $R_{\omega,\rho_T}$. Furthermore, the semi-conjugacy $h$ can be chosen so
 that it is fibre-respecting (i.e.\ $\pi_1\circ h = \pi_1$) and all fibre maps
 $h_\theta$ are order-preserving circle maps.  \listend
\end{thm}
\proof 
\alphlist
\item Suppose $T \in \Thom$ is regular and $\rho = \rho_T = \frac{k}{q} \omega +
  \frac{l}{q}$ for some $k,l,q \in \N$. Let $F :
  \R^2 \ra \R^2$ be a lift of $T^q$ with $\rho_F = k\omega$. Note that the
  rotation number of $F$ on the base is $q\omega$, i.e.\ $F :
  (\hat{\theta},\hat{x}) \mapsto (\hat{\theta} +
  q\omega,F_{\hat{\theta}}(\hat{x})$. Now let $\gamhat_0(\thhat) :=
  \frac{k}{q}\thhat$ and define 
  \[    
    A := \overline{ \bigcup_{n \in \Z}
    F^n(\hat{\Gamma}_0)} \ .
  \]
   Obviously $A$ is $F$-invariant.  Further, $A$ is contained in the strip $K
   := \{ (\thhat,\xhat) \mid \xhat \in [\frac{k}{q}\thhat - C,\frac{k}{q}
   \thhat +C] \}$, where $C$ is a suitable constant as in Def.\
   \ref{def:boundedorbits}. For any $\thhat \in \R$ the set $A_{\thhat}$ is
   the intersection of a nested sequence of compact sets and therefore
   non-empty. Thus $A$ projects down to all of \R. W.l.o.g.\ we can assume $k
   > 2C+1$, otherwise we replace $k$ and $q$ by sufficiently large multiples
   $lk$ and $lq$. 

  Now let $\widehat{B} := A^{\mathit{fill}}$. Due to the periodicity of $F$ and the
  definition of $A$ and $\widehat{B}$ we have $\widehat{B}_{\thhat +q} = \widehat{B}_{\thhat} + k$. Thus
  $\widehat{B}$ projects down to a simple strip $B$ in $\kreis_q \times \kreis_k$. As $B$ is contained in the projection of $K$ and $k>2C+1$, it is easy to see that there
  exists a simple curve in the complement.  $B$ is an invariant simple strip
  with respect to the projection of $F$ to $\kreis_q \times \kreis_k$, which
  completes the proof via Lemma \ref{lem:liftedgraph} .

\item
Now suppose the two rotation numbers are rationally independent. In order to
 show that $T$ and $R_{\omega,\rho_T}$ are semi-conjugate we will proceed in
 two steps. First, we show that there exists a family of disjoint minimal
 strips $(B_r)_{r\in\R}$ in $\kreis \times \R$ with $\That(B_r) =
 B_{r+\rho_T}$, $B_{r+n} = B_r + (0,n)$, and such that $r \mapsto B_r$ is
 strictly order-preserving. Then we use this to construct a semi-conjugacy $H$
 between $\That$ and $\widehat{R}_{\omega,\rho}$, which projects down to a
 semi-conjugacy $h$ between the original systems.

\ \\
\textbf{Step 1:} \ 
 Let $\That : \kreis \times \R \ra \kreis \times \R$ 
be a lift of $T$ with $\rho_{\That} = \rho_T =: \rho$. As
$|\Thatth^n(\xhat) - \xhat -n\rho| \leq C \ \forall n \in \N ,\ 
(\theta,\xhat) \in \kreis \times \R$, the sets 
\begin{equation}
    A_r := \overline{ \bigcup_{n\in\Z}
    \That^n(\kreis \times \{ r - n\rho \}) } 
\end{equation}
are all contained in $\kreis \times [r-C,r+C]$, and they are compact as the
intersection of a nested sequence of compact sets. As a direct consequence of
the definition $A_{r+n\rho} = \That^n(A_r) \ \forall n \in \Z$, the
periodicity of the fibre maps \Thatth\ implies $A_{r+n} = A_r + (0,n)$.

Furthermore, it follows directly from the monotonicity of the fibre
maps \Thatth\  
that the $A_r$ are ordered, i.e.\ $A_r \preccurlyeq A_s$ whenever $r \leq
s$. However, this order need not necessarily be strict.  
Therefore let
\begin{equation}
   \label{eq:Br}
    B_r := \overline{(\Phi^+_{A_r})^{-}}^{\mathit{fill}} \ .
\end{equation}
This means that we take the highest minimal strip contained in $\filled{A}$, compare
Section \ref{Minimalsets} .  These sets inherit the ordering property of the
sets $A_r$, as all the actions in (\ref{eq:Br}) preserve the order of sets. We
now aim to show that they are all disjoint, such that the ordering is
strict (i.e.\ $B_r \prec B_s$ whenever $r < s$).

\ \\
 To that end we claim that whenever $B_r \cap B_s \neq \emptyset$ for
 some $r<s$ there exists a regular invariant graph  
for $T$, contradicting the rational independence of $\omega$ and $\rho_T$:

Suppose $\exists \theta \in \kreis: (B_r)_\theta \cap (B_s)_\theta
\neq \emptyset$ and choose $r',s' \in \R$ such that   
$r+\delta \leq r' < s' \leq s-\delta$ for some $\delta > 0$. If we let 
\[
    \phiplus_r := \phiplus_{B_r} \ \mathrm{ and } \ \phimin_r :=
    \phimin_{B_r} \ \ \forall r \in \R 
\]
then $(B_r)_\theta \cap (B_s)_\theta \neq \emptyset$ implies
$\phiplus_r(\theta) \geq \phimin_s(\theta)$. 
Now, for any $n\in\Z$ with $n\rho \bmod 1 \in [-\delta,\delta] \bmod 1$
we have  
$B_{r-k} \leq B_{r'-n\rho} \leq B_{s'-n\rho} \leq B_{s-k}$ for some $k
\in \Z$. This means 
\[
    \phiplus_{r'-n\rho}(\theta) \geq \phiplus_{r-k}(\theta) =
    \phiplus_r(\theta) - k \geq \phimin_s(\theta ) - k =  
    \phimin_{s-k}(\theta) \geq \phimin_{s'-n\rho}(\theta)
\]
and thus $(B_{r'-n\rho})_\theta \cap (B_{s'-n\rho})_\theta \neq
\emptyset$. As $\That^n(B_{r'-n\rho}) = B_{r'}$ 
and $\That^n(B_{s'-n\rho}) = B_{s'}$ this means
$(B_{r'})_{\theta+n\omega} \cap (B_{s'})_{\theta+n\omega} \neq
\emptyset$.  

Consequently $\pi_1(B_{r'} \cap B_{s'})$ contains the set 
$\{ \theta + n\omega \mid n\rho \bmod 1 \in [-\delta,\delta] \}$, which 
is dense in \kreis\ as it is the projection of
$\{ R_{\omega,\rho_T}^n(0,0) \mid n \in \Z \} \cap (\kreis \times
[-\delta,\delta])$ onto 
the circle and $R_{\omega,\rho_T}$ is minimal. By the minimality of
the two strips $B_{r'}$ and $B_{s'}$  
this implies $B_{r'} = B_{s'}$.
Now, as we can take $r'=n\rho \bmod 1$ and $s'=k\rho \bmod 1$ for some
$n,k\in\Z$, this yields 
 $\That^{n-k}(B_{r'}) = B_{r'} + l$ for some $l \in \Z$. Therefore $B_{r'}$
projects down to an invariant strip for $T$. This completes Step 1. 

\ \\ \textbf{Step 2:} \ Now define $H : \kreis \times \R \ra \kreis \times
\R$ by $H(\theta,\xhat) := (\theta,H_\theta(\xhat))$ and
\[
    H_\theta(\xhat) \ := \ \sup\{r \in \R \mid \phiplus_r(\theta) \leq
    \xhat \} \ . 
\]
From the properties of the sets $B_r$ and their bounding graphs
$\phiplus_r$ it follows easily that  
$H \circ \That = \widehat{R}_{\omega,\rho} \circ H$, all $H_\theta$
are order-preserving (i.e.\ monotone) and 
that $H_\theta(\xhat + n) = H_\theta(\xhat) + n \ \forall n \in \Z, \
\xhat \in \R$. It remains to show that 
$H$ is continuous, i.e.\ that preimages of open sets are open. 

To that end consider an open set $U \subseteq \kreis \times \R$ and
$(\theta_0,\yhat) \in U$. There is a closed box $W :=
\overline{B_\delta(\theta_0)} \times [r,s]$ with $(\theta_0,\yhat) \in
W \subseteq U$. But this means $H^{-1}(U)$ contains the set $\{
(\theta,\xhat) \mid \theta \in B_\delta(\theta_0), \phiplus_r(\theta)
< \xhat < \phimin_s(\theta) \}$, which is an open neighbourhood of
$H^{-1}(\theta_0,\yhat)$ (openness following from the semi-continuity
of the graphs $\phiplus_r$ and $\phimin_s$).

This means $H$ is the required semi-conjugacy between the lifts of $T$
and $R_{\omega,\rho}$, and as  
$H_\theta (\xhat + n) = H_\theta(\xhat) + n$ we can project it down
to $\ntorus$ to obtain a semi-conjugacy $h$ between  
the original systems. 
\listend

\qed

\begin{bem}
  \alphlist
\item In case (a) it is also possible to show that the rotation
  numbers determine the numbers $p$ and $q$ in Def.\ \ref{def:tubes},
  as well as some further combinatorical properties of the invariant
  strips. A detailed description of this can be found in
  \cite{jaeger/keller:2003}.
\item In case (b) the semi-conjugacy $h$ is uniquely determined up to
  constant rotation in the direction of the fibres. This can be seen
  as follows: Suppose $h$ and $\tilde{h}$ have the properties
  mentioned in the theorem. Then first of all, as both must map the
  unique $T$-invariant measure $\mu$ to the Lebesque measure on
  $\ntorus$, we have $h_\theta = \tilde{h}_\theta + c_\theta$
  $m$-a.s., where the constant $c_\theta$ may a priori depend on
  $\theta$. But as $c_{\thom} = h_{\thom}(\Tth(x)) -
  \tilde{h}_{\thom}(\Tth(x) = h_\theta(x) +\rho_T -\tilde{h}_\theta(x)
  - \rho_T = c_\theta$ and the rotation on the base is ergodic, $c$ is
  $m$-a.s.\ constant. By continuity this extends to all fibres.
\listend
\end{bem}
Now we turn to the existence of wandering open sets. In dimension one,
the classical examples of circle homeomorphisms with irrational
rotation number and wandering open sets were given by Denjoy (see
\cite{demelo/vanstrien:1993}). They are constructed by ``blowing up''
one or more orbits of the respective irrational rotation. The
following lemma shows that this is also the only possibility of
producing such examples in the quasiperiodically forced case, with
orbits of points being replaced by orbits of line segments. As we will
see that irregular systems are always topologically transitive and can 
therefore not have any wandering open sets, it
suffices to consider the regular case.
\begin{lem}
  \label{lem:wanderingsets}
Suppose $T \in \Thom$ is semi-conjugate to the irrational torus translation
$R_{\omega,\rho}$ by a semi-conjugacy $h$ as in Thm.\ \ref{thm:poincare}
and let $U$ be a connected wandering open set. Then $h(U)$ is a straight line
segment. 
\end{lem}
\proof \\ Let $I \ssq \pi(U)$ and $\psi : I \ra \kreis$ be any continuous curve
contained in $U$. By $\Psi$ we denote the point set corresponding to
$\psi$. Now as $h$ is only a semi-conjugacy and therefore not neccessarily one
to one, the curves $R_{\omega,\rho:T}^n\circ h(\Psi)\ (n\in\N)$ may not be
pairwise disjoint. However, as $h$ preserves the order on each fibre they may
touch, but not cross each other. Obviously, the only curves in $\ntorus$ which
never cross any of their images under an irrational translation are straight
line segments. Thus any such curve contained in $U$ is mapped to a line
segment. Suppose that there are two curves $\psi_1,\psi_2$ in $U$ which are
mapped to different line segments. Then there exists a third curve $\varphi$
in $U$ which coincides with $\psi_1$ on one interval and with $\psi_2$ on
another. But then $h(\varphi)$ cannot be contained in a line segment, leading
to a contradiction. Therefore any curve, and in fact any point contained in $U$
must be mapped into the same line segment by $h$. 

\qed

\ \\ We close this section with a few remarks concering the
generalisation of Denjoy's Theorem (e.g.
\cite{demelo/vanstrien:1993,katok/hasselblatt:1997}) to
quasiperiodically forced systems. The existence of wandering open sets
can be excluded in a similar fashion to that for circle homeomorphims
by requiring that the map is sufficiently smooth. The condition needed
here is
\begin{equation}
  V(T) \ := \ \int_{\kreis} \textrm{Var} DT_\theta \ d\theta \ < \ \infty \ ,
\end{equation}
where $DT_\theta$ denotes the derivative of the fibre maps and
$\textrm{Var}(f)$ is the variation of the function $f$. However, unlike 
the one-dimensional case, excluding the possibility of wandering sets is not
enough to ensure the existence of a conjugacy to the irrational
translation. Therefore the following Denjoy-like statement, which is contained
in \cite{jaeger/keller:2003}, is somewhat weaker in nature:
\begin{thm}
  Let $T \in \Thom$ with $V(T) < \infty$ and suppose there exists no
  $T$-invariant strip. Then $T$ is topologically transitive. 
\end{thm}
An important step in the proof given in \cite{jaeger/keller:2003} is
to analyze the combinatorial behaviour of wandering open sets.
However, as irregular systems are always transitive (see next section)
and it therefore suffices to consider just the regular case, Lemma
\ref{lem:wanderingsets} can be used to simplify this part of the proof
significantly.


\section{The Irregular Case}
\label{Irregularcase}

First of all, we should mention that the irregular case is non-empty:
The simplest examples, which are due to Furstenberg, can be given by skew
translation of the torus, i.e. maps of the form
\begin{equation}\label{eq:skewtranslation}
  A : \ntorus \ra \ntorus \ , \ \ \ \thx \mapsto (\thom,x+a(\theta))
\end{equation}
were $a : \kreis \ra \kreis$ is continuous and homotopic to a constant. Now
there exist continuous (if $\omega$ is liouvillean even smooth) functions $\hat{a} : \kreis \ra \R$ with $\int_{\kreis}
\hat{a}(\theta) \ d\theta = 0$, such that the cohomological equation
\begin{equation}
  \hat{a}(\theta) \ = \ \phihat(\thom) - \phihat(\theta)
\end{equation}
has measurable, but no continuous solution (see \cite{mane:1987},
\cite{katok/robinson:2001}). This implies that the corresponding map $A$ given
by (\ref{eq:skewtranslation}) with $a = \pi \circ \hat{a}$ is minimal (see
\cite{katok/hasselblatt:1997}, Prop.\ 4.2.6). Thus there cannot be an
invariant strip. On the other hand, as $\varphi := \pi\circ\phihat$ gives an
invariant graph, $T$ can also not be semi-conjugated to an irrational
translation. Consequently, such a map must be irregular.

A slight modification of this construction also gives an example of an irregular
system without invariant graphs (Case II B below): If $\hat{a}$ is
replaced by $\hat{a} + \rho$ with some $\rho \in [0,1] \setminus \Q$,
then it is easily checked that $h\thx := (\theta,x-\varphi(\theta))$
satisfies $h^{-1}\circ R_{\omega,\rho} \circ h = T$ and therefore
gives an isomorphism between $T$ and the irrational translation. Thus
$T$ is uniquely ergodic with respect to the Lebesgue measure on
$\ntorus$ and therefore cannot have any invariant graphs.

The second interesting class of examples comes from the study of
matrix cocycles, given as maps
\begin{equation}
   T : \kreis \times \R^2 \ra \kreis \times \R^2 \ , \ \ \ (\theta,v) \mapsto (\thom,M(\theta)v)
\end{equation}
where $M : \kreis \ra \textit{SL}(2,\R)$ is continuous and homotopic
to $\textstyle\left(\begin{array}{rr} 1 & 0 \\ 0 & 1 \end{array}
\right)$. By the action of $M$ on the projective space $\Proj(\R^2)$,
such cocycles can be identified with quasiperiodically forced
M\"obius-transformations. In \cite{herman:1983} Herman showed that for
each pair $(\omega,\rho) \in \ntorus, \ \omega \in \kreis \setminus
(\Q/\Z)$ there exists a matrix cocycle over the rotation by $\omega$
with fibrewise rotation number $\rho$ and a positive Lyapunov
exponent. The later ensures the existence of two measurable invariant
graphs, such that the system cannot be semi-conjugated to an
irrational translation. On the other hand, if $\omega$ and $\rho$ are
rationally independent this excludes the existence of an invariant
strip. Thus the system must be irregular. In contrast to
Furstenberg's, these examples are analytic even if the rotation number
on the base is diophantine. They also show that the existence of
measurable invariant graphs does not require the rational dependence
of the rotation numbers.

The following theorem gives at least some information about the
topological dynamics in the irregular case:

\begin{thm}
If $T \in \Thom$ is irregular then it is topologically transitive.
\end{thm}
\proof \\ First we need to present a number of notations and conventions
for continuous curves $\varphi : I \ra \kreis$ where $I = [a,b] \ssq
\kreis$ is a compact interval. The point set which is given by
the graph of such a curve will be denoted by the corresponding capital letter
(as for invariant graphs). If $\varphi : I \ra \kreis$ is a continuous curve
the set $T(\Phi)$ is the point set of a continuous curve as well, which we will
denote by $T\varphi : I + \omega \ra \kreis$. Given any point $\xhat \in
\piin(\varphi(a))$ and any lift $\hat{I} = [\hat{a},\hat{b}]$ of the interval
$I$, there exists a unique lift $\phihat: \hat{I} \ra \R$ of $\varphi$ which
satisfies $\phihat(\hat{a}) = \xhat$. Using any such lift we can define the
\textit{winding number} of the curve $\varphi$ by $k(\varphi) :=
\phihat(\hat{b}) - \phihat(\hat{a})$. Further let $v(\varphi) := \max_{\thhat
\in \hat{I}} \phihat(\thhat) - \min_{\thhat \in \hat{I}} \phihat(\thhat)$. It is
easy to see that if $\psi:I' \ra \kreis$ is another curve with $I' \ssq I$ and
$k(\psi) \geq v(\varphi) + 1$, then $\Psi$ must intersect $\Phi$.

We have to show that given any two open sets $U,V$ there exists some
$n\in\N$, such that $T^n(U) \cap V \neq \emptyset$. It suffices to consider
the case where $U$ and $V$ are open rectangles, i.e.\ $U = I_1\times J_1$ and
$V = I_2 \times J_2$ for some open intervals $I_1,I_2,J_1,J_2 \ssq \kreis$.
Now take a compact interval $I \ssq I_2$ and any continuous curve $\varphi : I
\ra \kreis$ with $\Phi \ssq V$. For some $n \in \N$, the interiors of the
intervals $I, T^{-1}(I) \ld T^{-n}(I)$ cover the whole circle. Consequenty
there is some $\delta > 0$ such that any interval $I'$ of length smaller
or equal to $\delta$ is fully contained in one of these intervals. By
Thm.\ \ref{thm:boundedness} there exists both an orbit which is
$\rho_T$-bounded above as well as one which is $\rho_T$-unbounded
above. Therefore we can choose an interval $I' = [\theta_1,\theta_2] \ssq I_1$
of length smaller or equal to $\delta$, such that the orbits on the
$\theta_1$-fibre are $\rho_T$-bounded above whereas the orbits on the
$\theta_2$-fibre are not. This means that if we take any continuous curve
$\psi : I' \ra \kreis$ with $\Psi \ssq U$, we can find some
$m \in \N$, such that $k(T^m\psi) \geq \max_{i=0}^n v(T^i\phi) +1$. As
$T^m(I')$ is contained in $T^{-j}(I)$ for some $j \in \{0 \ld n\}$, this
implies that $T^m(\Psi)$ intersects $T^{-j}(\Phi)$. But as $\Psi \ssq U$ and
$\Phi \ssq V$ we are finished.

\qed

\ \\
Combining the results from the last section with Thm.\ \ref{thm:furstenberg},
we obtain the following basic classification:

\begin{picture}(400,300)(30,0)
\put(10,10){\framebox(360,260)}

\put(90,10){\line(0,1){260}}
\put(10,230){\line(1,0){360}}
\put(10,120){\line(1,0){360}}
\put(230,10){\line(0,1){260}}

\put(140,240){\shortstack{regular \\ (I) }}
\put(280,240){\shortstack{irregular \\ (II) }}
\put(30,160){\shortstack{invariant \\ graphs \\ (A) }}
\put(30,45){\shortstack{no \\ invariant \\ graphs \\ (B) }}

\put(115,155){\shortstack{invariant strips \\ \ \\ \ \\ \ \\ rationally
    dependent \\ rotation numbers }} 
\put(110,40){\shortstack{semi-conjugacy
    to the \\ irrational translation \\ \ \\ \ \\ \ \\ rationally independent
    \\ rotation numbers }} 
\put(260,45){\shortstack{isomorphism to a \\ uniquely ergodic \\ skew
    translation \\ \ \\ \ \\ \ \\ transitive }}
 \put(275,175){\shortstack{transitive }}
\end{picture}

\ \\
Let us examine Case II(B) more closely: Suppose $T \in \Thom$ is
isomorphic to the uniquely ergodic skew translation $A$ via the
isomorphism $h$. In general, an isomorphism of two systems on the
torus cannot be lifted to $\kreis \times \R$ in any obvious way.
However, although not explicitly mentioned in \cite{furstenberg:1961},
the construction of $h$ in the proof of Thm.\ \ref{thm:furstenberg}
gives a lot of additional information about $h$. It is the projection
of a function $\hat{h} : \kreis \times \R \ra \kreis \times \R$, which
is an isomorphism between a lift $\That$ of $T$ and a skew translation
$\Ahat$ which projects down to $A$. Further, it is of the form
$\hat{h}\thx = (\theta,\hat{h}_\theta(x))$ with all fibre maps
$\hat{h}_\theta$ being monotone and continuous, and finally $\hat{h}$
leaves all integer lines $\kreis \times \{n\}$ invariant. This last
property ensures that $\hat{h}$ preserves the rotation number, in the
sense that $\rho_{\That} = \int_\kreis \hat{a}(\theta) \ d\theta$.
This leaves the following possibilities:
\romanlist  
  \item $\hat{a}$ is cohomologous to a constant, which is neccessarily equal
  to $\rho_{\That}$, i.e.\
  \begin{equation}\label{eq:cohomology}
        \hat{a}(\theta) \ = \ \phihat(\thom) - \phihat(\theta) + \rho_{\That}
  \end{equation}
  for some measurable function $\phihat : \kreis \ra \R$. As we have
  seen in the discussion of Furstenberg's examples above, this means
  that $T$ is isomorphic to $R_{\omega,\rho_T}$ by a fibre-respecting
  isomorphism $h$. As isomorphy preserves ergodicity and $h$ leaves
  $\lambda$ invariant, $R_{\omega,\rho_T}$ must be ergodic w.r.t.\ 
  Lebesgue measure as well, and therefore irrational. In particular,
  $\omega$ and $\rho_T$ cannot be rationally dependent in this case.
\item There exists no solution to the cohomological equation
  (\ref{eq:cohomology}). In this case it is hard to say anything further, see
  questions (a) and (b) below. 
\listend
We close with three questions related to this discussion, which we have to leave
open here:
\begin{questions}
  \alphlist
  \item
  Is is possible that a system $T \in \Thom$ is isomorphic to an irrational
  translation $R_{\omega,\rho}$ with $\rho$ not being contained in the module
  $\{ \rho_T + k \omega \mid k \in \Z \}$? (Note that two irrational
  translations $R_{\omega,\beta}$ and $R_{\omega,\beta+k\omega}$ from the same
  module are always isomorphic by $h(\theta,x) = (\theta,x+k\theta)$).
  \item
  Is it possible to have rationally related rotation numbers in Case II B?
  \item Is it possible to have transitive but not minimal dynamics in the
  irregular case?
  \listend
\end{questions}

\ \\
\textbf{Acknowledgements.} We would like to thank Sylvain Crovisier
for interesting discussions on the subject. Tobias Jaeger acknowledges
support from the DFG and also profited from a Marie Curie Fellowship
during a visit to the University of Surrey.

\bibliography{class} \bibliographystyle{alpha}

\begin{thebibliography}{SFGP02}

\bibitem[AK]{avila/krikorian}
A.~Avila and R.~Krikorian.
\newblock Reducibility or non-uniform hyperbolicity for quasiperiodic
  {S}chr{\"o}dinger cocycles.
\newblock Preprint.

\bibitem[DGO89]{ding/grebogi/ott:1989}
Mingzhou Ding, Celso Grebogi, and Edward Ott.
\newblock Evolution of attractors in quasiperiodically forced systems: from
  quasiperiodic to strange nonchaotic to chaotic.
\newblock {\em Physical Review A}, 39(5):2593--2598, 1989.

\bibitem[dMvS93]{demelo/vanstrien:1993}
W.~de~Melo and S.~van Strien.
\newblock {\em One-dimensional dynamics}.
\newblock Springer, 1993.

\bibitem[FKP95]{feudel/kurths/pikovsky:1995}
Ulrike Feudel, J{\"u}rgen Kurths, and Arkady~S. Pikovsky.
\newblock Strange nonchaotic attractor in a quasiperiodically forced circle
  map.
\newblock {\em Physica D}, 88:176--186, 1995.

\bibitem[Fur61]{furstenberg:1961}
H.~Furstenberg.
\newblock Strict ergodicity and transformation of the torus.
\newblock {\em American Journal of Mathematics}, 83:573--601, 1961.

\bibitem[Her83]{herman:1983}
Michael~R. Herman.
\newblock Une {m\'{e}thode} pour minorer les exposants de {{L}yapunov} et
  quelques exemples montrant le {caract\`{e}re} local d'un {th\'{e}or\`{e}me}
  {d'Arnold} et de {Moser} sur le tore de dimension 2.
\newblock {\em Commentarii Mathematici Helvetici}, 58:453--502, 1983.

\bibitem[Her86]{herman:1986}
M.~Herman.
\newblock Construction of some curious diffeomorphisms of the {R}ieman sphere.
\newblock {\em Journal of the London Mathematical Society}, 34:375--384, 1986.

\bibitem[JK]{jaeger/keller:2003}
T.~J{\"a}ger and G.~Keller.
\newblock The denjoy type-of argument for quasiperiodically forced circle
  diffeomorphisms.
\newblock Preprint.

\bibitem[KH97]{katok/hasselblatt:1997}
A.~Katok and B.~Hasselblatt.
\newblock {\em Introduction to the Modern Theory of Dynamical Systems}.
\newblock Cambridge University Press, 1997.

\bibitem[KR01]{katok/robinson:2001}
A.~Katok and E.~A. Robinson.
\newblock Cocycles, cohomology and combinatorial constructions in ergodic
  theory.
\newblock {\em Proceedings of Symposia in Pure Mathematics}, 69:107--173, 2001.

\bibitem[Ma{\~n}87]{mane:1987}
Ricardo Ma{\~n}{\'e}.
\newblock {\em Ergodic Theory and Differentiable Dynamics}.
\newblock Springer, 1987.

\bibitem[MZ89]{misiurewicz/ziemian:1989}
M.~Misiurewicz and K.~Ziemian.
\newblock Rotation sets for maps of tori.
\newblock {\em Journal of the London Mathematical Society}, 40(3):490--506,
  1989.

\bibitem[RT86]{rhodes86}
F.~Rhodes and C.L. Thompson.
\newblock Rotation numbers for monotone functions of the circle.
\newblock {\em J. Lond. Math. Soc.}, 34:360--368, 1986.

\bibitem[RT91]{rhodes91}
F.~Rhodes and C.L. Thompson.
\newblock Topologies and rotation numbers for families of monotone functions on
  the circle.
\newblock {\em J. Lond. Math. Soc.}, 43:156--170, 1991.

\bibitem[SFGP02]{stark/feudel/glendinning/pikovsky:2002}
J.~Stark, U.~Feudel, P.A. Glendinning, and A.~Pikovsky.
\newblock Rotation numbers for quasi-periodically forced monotone circle maps.
\newblock {\em Dynamical Systems}, 17(1):1--28, 2002.

\bibitem[Sta03]{stark:2003}
J.~Stark.
\newblock Transitive sets for quasiperiodically forced monotone maps.
\newblock {\em Dynamical Systems}, 18(4):351--364, 2003.

\end{thebibliography}

\end{document}